%% file: NonMonic.tex
\documentclass{amsart}
\usepackage{graphics}
\usepackage{amsfonts}
\usepackage{amssymb}
\usepackage[mathcal]{eucal}
\usepackage{mathrsfs}
\usepackage{latexsym}
\usepackage{amsthm}
\usepackage{amsmath,amscd}
\usepackage{enumerate}
\usepackage{epic}
\usepackage{hyperref}
\input{Macros}

\title[Nonmonic factoring]{An algorithm for factoring non-monic quadratic polynomials\\Or: How I learned to stop using the quadratic formula and love undoing FOIL}

\author[C.T. Bruns]{Corey~T. Bruns}
\address{Department of Mathematics\\University of Colorado\\CB 395\\Boulder, CO USA}
\email{brunsc@colorado.edu}
\urladdr{\url{http://math.colorado.edu/~brunsc}}

\thanks{\copyright 2009 Corey T.~Bruns\\
This work is licensed under the Creative Commons Attribution-Noncommercial-Share Alike 3.0 United States License. To view a copy of this license, visit   \url{http://creativecommons.org/licenses/by-nc-sa/3.0/us/} or send a letter to Creative Commons, 171 Second Street, Suite 300, San Francisco, California, 94105, US}

\begin{document}
\maketitle
\begin{abstract}
{We give an algorithm for factoring quadratic polynomials over any UFD, $\Z$ in particular.  We prove the correctness of this algorithm and give examples over $\Z$ and $\Z\left[i\right]$.}
\end{abstract}

\input{nonmonicfactoring}

\bibliographystyle{amsalpha}
\bibliography{biblio}

\end{document}

%% file: Macros.tex
\usepackage{ucs}
\usepackage[utf8x]{inputenc}
\long\def\symbolfootnote[#1]#2{\begingroup\def\thefootnote{\fnsymbol{footnote}}\footnote[#1]{#2}\endgroup}

\newcommand{\Z}{\mathbb{Z}}

\newtheorem{theorem}{Theorem}[section]

\theoremstyle{definition}

\newcommand{\sm}{\smallsetminus}

\renewcommand{\phi}{\varphi}

\newcommand{\pa}[1]{\left(#1\right)}

\newcommand{\s}[1]{\left\{#1\right\}}

\renewcommand{\epsilon}{\varepsilon}

\newcommand{\br}[1]{\left[#1\right]}

%% file: nonmonicfactoring.tex
\section{Introduction}
Factoring quadratic trinomials (over $\Z$) (e.g. starting with $4x^2+8x+3$ and ending up with $\pa{2x+1}\pa{2x+3}$) is a major topic in elementary algebra; however, most treatments of the subject don't treat the non-monic case \footnote{Recall that a monic polynomial is one in which the leading coefficient is 1.} in a systematic way.  Some books advise trying all possible factorizations of the leading term and all possible factorizations of the constant term as the terms of the binomial factors, and seeing which one works by multiplying out the factors.  Even in the simple example above, this leads to checking three possibilities, which is not too bad; on average this means performing FOIL \footnote{First, Outside, Inside, Last; a mnemonic device for the result of the three applications of the distributive property in the expansion of a  product of two binomials} 1.5 times. If the coefficients are not prime powers, the number of possibilities is much larger; something like $6x^2+7x-24$ gives sixteen possibilities, which is enough to make even the most ardent student think there must be a better way.  Often this better way is to use the quadratic formula and work backwards from the roots to the factorization; indeed, as an algorithm, this is probably optimal.  However, it's probably better to have done some factoring of these types of polynomials before the quadratic formula is introduced.  As well, since we are really trying to factor a polynomial in $\Z\br x$, it's good to do this without resorting to square roots.  The method we present does this, and has the additional benefit of promoting a better understanding of the process of factorization in general, even the monic case. Furthermore, this algorithm works for polynomials over any unique factorization domain; one must find greatest common divisors (GCD) in the original ring, so a Euclidean domain of coefficients makes this a more reasonable task.

I don't claim to have invented this technique; another example is at \\\url{http://mathforum.org/library/drmath/view/52878.html}.  My purpose is to give a proof of the correctness of this technique and to raise its profile so that it may make it into algebra texts.  

\section{The Theory}\label{sec:theory}
Our goal is to begin with a trinomial and write it as the product of two binomials, if possible.  The reverse process, which students have doubtlessly done many times by this point in their careers, is to begin with a product of two binomials, and by performing FOIL and combining like terms, typically the Outside and Inside products.
We use the phrase ``undoing FOIL'' as it expresses the intent even though there are additional steps; we must first undo the combining like terms step.  Indeed, this step is the meat of the method. Of course there are many ways that the middle term of the trinomial can be written as a sum; the trick here is to do it the right way, that is, the way that allows us to actually undo the FOIL.
 We will introduce some conventions here; we assume that the trinomial is of the form $ax^2+bx+c$, though this applies to other variables and similar types of trinomial, such as $ax^2+bxy+cy^2$ or $ax^4+bx^2t^3+ct^6$. It is critical that $a\neq 0$ and $\mbox{GCD}\pa{a,b,c}=1$, and interesting if $a\neq 1$ and $bc\neq 0$.  (For general UFDs, substitute your symbols of choice for the additive and multiplicative identities.)  The final factorization will be symbolized $\pa{Ax+B}\pa{Cx+D}$,  where $\mbox{GCD}\pa{A,B}=\mbox{GCD}\pa{C,D}=1$; mutatis mutandis if you have a different form in mind.

\begin{theorem} Let $R$ be a UFD (e.g. $\Z$).  Let $a\in R\sm\s{0,1}$; $b,c\in R\sm\s 0$ be such that $\mbox{GCD}\pa{a,b,c}=1$. Then the polynomial $ax^2+bx+c\in R\br x$  factors into two linear factors if and only if there exist $b_1, b_2\in R$ such that $b_1 b_2=ac$ and $b_1 + b_2=b$. Furthermore, let $A=\mbox{GCD}\pa{a,b_1}$, $B=\mbox{GCD}\pa{b_2,c}$, $C=\mbox{GCD}\pa{a,b_2}$, and $D=\mbox{GCD}\pa{b_1,c}$; then $ax^2+bx+c=\pa{Ax+B}\pa{Cx+D}$.
\end{theorem}
\begin{proof}
So, taking   $\pa{Ax+B}\pa{Cx+D}$, and applying FOIL to it, we get $$ACx^2 +BCx +ADx+BD=ACx^2+\pa{BC+AD}x + BD.$$  So, in theory, all we have to do is solve the following system of equations: $$\left\{\begin{array}{rrr}AC&=&a\\ BC+AD&=&b\\ BD&=&c\end{array}\right.$$ Clearly, $ax^2+bx+c$ factors into two linear factors if and only if this system has a solution. We do this by combining the first and third equations here; getting $ABCD=ac$; we regroup this as $\pa{BC}\pa{AD}=ac$; then change variables $b_1=AD$ and $b_2=BC$; then we just have $$\left\{\begin{array}{rrr}b_1 b_2&=&ac\\b_1 + b_2&=&b\end{array}\right.$$  This is a good simplification; equivalent to factoring $x^2+ bx + ac$, so we are back to the monic case.  More explicitly, we just need to figure out which pair of factors of $ac$ adds up to $b$.\footnote{This can be systematically done by our hypothetical student by exhaustively listing the factor pairs of $ac$ and summing each pair to compare to $b$. Note that each factor pair will have opposite signs if $ac<0$, so that this summing may be better explained as subtracting in this case. In the more advanced case, if $R$ is a Euclidean Domain, this listing of factor pairs in order (by the smaller of the norms of the factors) can still be accomplished.  For other $R$, one could factor these numbers in $R$ and partition the factorization systematically to find $b_1$ and $b_2$.}  Now that we have $b_1$ and $b_2$, we need to get back to $A, B, C, \mbox{and } D$.\footnote{This is the part where we hope that the students may learn something about factoring.}  The change of variable was chosen so that we may now take our original trinomial and write it as $ax^2+b_1 x + b_2 x + c$; that is, the transformed system is what we really wanted in order to undo the combination of like terms!  At this point, we can proceed to factor by grouping; since $a=AC$ and $b_1=AD$, the GCD of $ax^2+b_1 x$ is $Ax$; similarly the GCD of $b_2 x + c$ is $B$, etc.
So we have the following steps:  $$\begin{array}{ll}&ax^2+bx+c  \\=&ax^2+b_1x+b_2x+c\\=&Ax(Cx+D)+B(Cx+D)\\=&(Ax+B)(Cx+D)\end{array}$$
\end{proof}

\section{The Practice}
In practice, it may not be the case that $\mbox{GCD}\pa{a,b,c}=1$; in this case, the first factorization step is to factor out that GCD; we then proceed on the resulting trinomial.  (If this is skipped, then the formulas for $A, B, C, \mbox{and } D$ at the end of the previous section will fail.)  So, when presenting this method, it is important to emphasize that this is, as always, the first step.
The procedure can be explained with something like the following:

\begin{quote}We are trying to ``undo'' a FOIL operation.  Start by multiplying the leading coefficient and the constant term.  If this number is positive, we need to find a pair of its factors that add up to the middle term's coefficient; if it's negative the factors need to subtract to the middle term's coefficient.  These factors are how we'll split up the middle term. \end{quote}

At this point, show an example.  Now that we have four terms, we can label each one as the F, O, I, and L term and proceed from here, either by grouping, or by the graphical means shown in the following section as the instructor prefers.
\section{Elementary Examples}

We give a couple of concrete examples with varying signs here.  We will redraw figures  as each additional mark is made; during a class each addition should be explained and just added to the existing figure. 

For our first example, with all signs the same; we revisit $4x^2+8x+3$; $a=4$, $b=8$, and $c=3$.  $ac=4\cdot 3=12$, so we list the factorizations of $12$; $1\cdot 12$, $2\cdot 6$, and $3\cdot 4$.  Since $ac>0$, we need to add these factors; looking for a sum of $b=8$.  The second pair, $2\cdot 6$ does so; so we may rewrite the original polynomial as $4x^2+2x+6x+3$. This time we will proceed by grouping; $4x^2+2x+6x+3=2x\pa{2x +1} + 3\pa{2x +1}=\pa{2x +3}\pa{2x+1}$.

We note that the procedure for factoring $4x^2-8x+3$ is almost identical; since $ac>0$, we still look for its factors that add up to $8$.  Technically, we ought to be looking for a pair of negative factors that add up to $-8$, but it may be easier to deal with the signs later.  The difference is that the rewritten polynomial is $4x^2-2x-6x+3$ and the resulting factorization is $\pa{2x-3}\pa{2x-1}$.

For a second example, with $ac<0$, we factor $6x^2+7x-24$.  Here $ac=-144$, so that the pair of factors will have opposite signs.  In practical terms, this means that we will look at positive factors of $144$ and look for a difference of $7$.  The factorization of $144$ that has this property is $16\cdot 9$.  Since the two factors of $-144$ should add to $7$, we actually use $16\cdot -9$, resulting in a rewritten polynomial of $6x^2+16x-9x-24$.  This time we will use a graphical representation of the remaining process.  We write the four terms in the cells of a $2\times 2$ array as follows:

$$\begin{array}{|c|c|}\hline
6x^2 & 16x\\\hline
-9x & -24\\\hline\end{array}$$

The purpose of this array is that the GCD of the rows and columns gives us part of the factorization; we will write these above the columns and to the left of the rows.  We can also choose the signs here; each term in the box ought to be the product of the GCDs above and to the left of it; indeed, after finding one GCD, the rest may be found with division.  We'll start with the GCD of the first column, which is $3x$, and write it as follows:

$$\begin{array}{|c|c|}

3x& \\\hline
6x^2 & 16x\\\hline
-9x & -24\\\hline\end{array}$$

Then we calculate $6x^2/3x=2x$ and write that to the left of the first row:

$$\begin{array}{c|c|c|}

& 3x& \\\hline
2x &6x^2 & 16x\\\hline
&-9x & -24\\\hline\end{array}$$

Then we calculate $16x/2x=8$ and $-9x/3x=-3$:

$$\begin{array}{c|c|c|}

& 3x&8 \\\hline
2x &6x^2 & 16x\\\hline
-3&-9x & -24\\\hline\end{array}$$

And now we can read off the factors: $\pa{3x+8}\pa{2x-3}$.

For our last example in $\Z[x]$, we will take the polynomial $15x^2-29x-14$.  Again $ac<0$, in particular it is $-210$.  So we need a pair of factors of $210$ that have a difference of $29$.  Systematically we proceed as follows:\\  $210-1=209$\\ $105-2=103$\\ $70-3=67$\\ $42-5=37$\\ $35-6=29$\\  There are of course other factor pairs, but we stop here as we've found the pair that will work.  With lucky guesses, it may take even fewer guesses; in this case the first two could reasonably be skipped.  In this case, the two $x$ terms must add up to $-29x$, so we take the factorization $-35\cdot 6$, resulting in $15x^2-35x+6x-14$.  We put this into the array form:

$$\begin{array}{|c|c|}
\hline
15x^2 & -35x\\\hline
6x & -14\\\hline
\end{array}$$

We find the GCD of the first row this time (we may choose any row or column to begin).
$$\begin{array}{c|c|c|}
\hline
5x&15x^2 & -35x\\\hline
&6x & -14\\\hline
\end{array}$$

Now calculate $15x^2/5x=3x$:
$$\begin{array}{c|c|c|}
&3x &\\ \hline
5x&15x^2 & -35x\\\hline
&6x & -14\\\hline
\end{array}$$
And finally $-35x/5x=-7$ and $6x/3x=2$:

$$\begin{array}{c|c|c|}
&3x &-7\\ \hline
5x&15x^2 & -35x\\\hline
2&6x & -14\\\hline
\end{array}$$

And we have our factorization: $\pa{3x-7}\pa{5x+2}$.

\section{Using some other UFD}

The theory is no different than in section \ref{sec:theory}.  If the ring is not ordered (as in the ensuing example), care is required in making sure that all factorizations of $ac$ are tested.  Since we intend to come up with a factorization into irreducible polynomials, this requires a UFD so that the factorization is unique.\footnote{Well, unique up to associates, which is the phrase that the U in UFD {\em really} stands for.}  Fortunately, all UFDs are also GCD domains, that is, any two elements have a unique (up to associates) greatest common divisor.\footnote{Here, {\it greatest} is under the divisibility ordering; that is $a$ is the GCD of set $X$ if it is a divisor of every element of $X$ and any other $b$ that divides every element of $X$ also divides $a$.}   It's also convenient that $R\br x$ is a UFD if and only if $R$ is a UFD (Chapter 8 of \cite{DandF} for both of these facts).  

We give an example in $\Z\br i \br x$.  We start with the polynomial 
$$\pa{2+4i}x^2 + \pa{7+5i}x +10.$$  We first notice that the coefficients have a common factor, $1+i$, so we factor that out first, obtaining
$$\pa{1+i}\br{ \pa{3+i}x^2 + \pa{6-i}x+\pa{5-5i} }.$$
So now we only deal with $\pa{3+i}x^2 + \pa{6-i}x+\pa{5-5i}$.  Now $ac=20-10i$; so we need a pair of factors that add up to $6-i$.  We factor $$20-10i=i\pa{1+i}^2\pa{1+2i}^2\pa{2+i},$$ where $i$ is a unit and the rest are irreducible.
  The fact that $\Z\br i$ is not ordered makes the next step somewhat less systematic; a different order of writing the factorization above will lead to considering factor pairs in a different order.
 We search for a factor pair with a sum of (an associate of) $6-i$ by partitioning this set; we get that $i\pa{1+2i}^2 +\pa{1+i}^2\pa{2+i}=\pa{-4-3i}+\pa{-2+4i}=-6+i=-1\pa{6-i}$, so we get $4+3i$ and $2-4i$ as the factors that add up to $6-i$.  Then we can rewrite our polynomial as 
$$\pa{3+i}x^2 + \pa{4+3i}x+\pa{2-4i}x+\pa{5-5i}.$$  If we proceed by grouping, we get $$\pa{1+2i}x\br{\pa{1-i}x+\pa{2-i}} + \pa{1+3i}\br{\pa{-1-i}x+\pa{-1-2i}}.$$
Then the terms in brackets are associates of each other; we factor out a $-i$ from the latter to get:
$$\pa{1+2i}x\br{\pa{1-i}x+\pa{2-i}} + \pa{3-i}\br{\pa{1-i}x+\pa{2-i}}.$$
And then finally, restoring the $1+i$ factor:
$$\pa{1+i}\br{\pa{1+2i}x+\pa{3-i}}\br{\pa{1-i}x+\pa{2-i}}.$$